\RequirePackage{fix-cm}
\documentclass[smallextended]{svjour3}       
\smartqed  
\usepackage{amsmath}
\usepackage{amsfonts}
\newcommand{\C}{\mathbb {C}}
\begin{document}

\title{Special classes of $q$-bracket operators
}


\author{Tanay Wakhare}


\institute{Tanay Wakhare \at
              University of Maryland, College Park, MD 20742 \\
              Tel.: +1-240-224-3900\\
              \email{twakhare@gmail.com} \\
	ORCID ID 0000-0003-2397-4185
}

\date{Received: date / Accepted: date}

\maketitle

\begin{abstract}
We study the $q$-bracket operator of Bloch and Okounkov, recently examined by Zagier and other authors, when applied to functions defined by two classes of sums over the parts of an integer partition. We derive convolution identities for these functions and link both classes of $q$-brackets through divisor sums. As a result, we generalize Euler's classic convolution identity for the partition function and obtain an analogous identity for the totient function. As corollaries, we generalize Stanley's Theorem on the number of ones in all partitions of $n$, and provide several new combinatorial results.
\keywords{partitions \and $q$-series \and multiplicative number theory \and additive number theory}
 \subclass{05A17 \and 11P81 \and 11N99}
\end{abstract}

\section{Introduction}
\label{Introduction}

Let $\mathcal P$ denote the set of integer partitions $\lambda=(\lambda_1, \lambda_2,\ldots, \lambda_r)$ with $\lambda_1\geq \lambda_2\geq \ldots \geq \lambda_r$, including the empty partition $\emptyset$. The $q$-bracket operator $\left< f \right>_q$ of a function $f: \mathcal P \to \mathbb C$ was introduced by Bloch and Okounkov in 2000 [3], and represents the average value of $f$ over all partitions (in the sense of statistical physics).  The $q$-bracket is defined in \cite[(Definition 1)]{Schneider1} by
\begin{equation}\label{1.1}
\langle f\rangle_q := \frac{\sum_{\lambda \in \mathcal{P}} f(\lambda)q^{|\lambda|}}{\sum_{\lambda \in \mathcal{P}}q^{|\lambda|}} = (q;q)_{\infty} \sum_{n=0}^{\infty}q^n \sum_{\lambda \vdash n}f(\lambda)\in \C[[q]],
\end{equation}
in which $\lambda \vdash n$ indicates that $\lambda \in \mathcal P$ is a partition of $n\geq 0$, and $\lambda_i \in \lambda$ means $\lambda_i$ is a part of $\lambda$. The $q$-bracket has been studied by Zagier \cite{Zagier1}, Griffin-Jameson-Trebat-Leder \cite{Griffin}, and Schneider \cite{Schneider1} due to connections to quasimodular forms, $p$-adic modular forms and phenomena in partition theory.

We will denote by $\sum_{\lambda \vdash n }f(\lambda)$ the sum of $f$ over all partitions of $n$. Moreover, for an arithmetic function $a: \mathbb Z_{\geq 0} \to \mathbb C$, in this paper we usually take $f(\lambda)$ to be of the form either $\sum_{\lambda_i \in \lambda} a(\lambda_i)$, with the sum taken over the parts $\lambda_i$ of $\lambda$, or $\sum_{\substack{\lambda_i \in \lambda\\ \lambda_i\  \text{distinct}}} a(\lambda_i)$, where the sum is over the distinct parts of $\lambda$ (i.e. the parts without repetition). For instance, $\sum_{\substack{\lambda_i \in (1,1,2,4,4,5) \\ \lambda_i \text{distinct} }}a(\lambda_i) = a(1)+a(2)+a(4)+a(5)$. As usual, we write $\sum_{d|n}f(d)$ to denote the sum of $f$ over all positive divisors of $n\geq 1$, including $1$ and $n$. We let $p(n)$ denote the number of unrestricted integer partitions of $n\geq 0$, with the convention $p(0):=1$.

\section{Sums over all parts}
\label{Sums1}
When our function $f$ splits over every part of a partition, its $q$-bracket turns out to have a very nice representation.
\begin{theorem}\label{Theorem1}
For $\lambda \in \mathcal P$ let $f(\lambda)=\sum_{\lambda_i \in \lambda} a(\lambda_i)$. Then, using the notation introduced above, we have
$$
\left< f \right>_q=\sum_{n=1}^{\infty} \frac{a(n) q^n}{1-q^n}.
$$
\end{theorem}
\begin{proof}
We first find a product that will generate $\sum_{n=0}^{\infty}q^n \sum_{\lambda \vdash n}f(\lambda)$, where $f(\lambda)=\sum_{\lambda_i \in \lambda}a(\lambda_i)$ for convenience. We begin with the usual \cite[pp. 3-4]{Andrews2}
\begin{equation}
\sum_{n=0}^{\infty}p(n)q^n = \frac{1}{(q;q)_\infty} = (1+q+q^2+\ldots)(1+q^2+q^4+\ldots)\cdots.
\end{equation}
We introduce a new variable $t$ and study 
\begin{align*}
\phi(t,q) &:= (1+t^{a(1)} q+t^{2a(1)}q^2+\ldots)(1+t^{a(2)} q^2+t^{2a(2)}q^4+\ldots) \times\cdots \\
&\times (1+t^{a(n)} q^n+t^{2a(n)}q^{2n}+\ldots)\times \cdots.
\end{align*}
Now every coefficient of $q^n$ will be a polynomial in $t$, with exponent $m$ in $t^m$ equal to $a(n_1)+a(n_2)+\ldots$ for $n_1, n_2, \ldots, n_i$ that form a partition of $n$. Therefore, by taking the partial derivative of $\phi(t,q)$ with respect to $t$ and then setting $t=1$ we will sum the powers of each $t$-polynomial and make that the new coefficient of $q^n$. We then have the relation 
\begin{equation}
\frac{\partial \phi(t,q)}{\partial t}\Bigr|_{t=1} = \sum_{n=0}^{\infty}q^n \sum_{\lambda \vdash n}f(\lambda) = \frac{1}{(q;q)_{\infty}}\langle f\rangle_q.
\end{equation}
Evaluating the partial by summing each term of $\phi$ as a geometric series, and taking the logarithmic derivative of  
\begin{equation}
\phi(t,q) = \prod_{n=1}^{\infty}\left(\frac{1}{1-t^{a(n)}q^n}\right)
\end{equation}
while noting $\phi(1,q) = \frac{1}{(q;q)_\infty}$, yields the desired theorem.
\end{proof}

We note that the case $a(n)=n^\alpha$ was treated in \cite{Han1}. We also note that the expression on the right-hand side of Theorem \ref{Theorem1} is simply a Lambert series \cite[(27.7.5)]{NIST:DLMF}. The case $a(\lambda_i)=\lambda_i^{2k-1}$ was considered by Zagier \cite[(42)]{Zagier1} as his ``moment function". We find that the resulting Lambert series is equal to $\sum_{n=0}^{\infty}\sigma_{2k-1}(n)q^n$, which arises as part of the Fourier expansion of weight $2k$ Eisenstein series. This shows the connection of certain $q$-brackets with quasimodularity, and highlights the special nature of the moment function: it is the only function additive over parts of partitions whose $q$-bracket will generate weighted Eisenstein series.

\begin{theorem}\label{Theorem2}
We have the convolution identity
\begin{equation}
\sum_{\lambda \vdash n}\sum_{\lambda_i \in \lambda}a(\lambda_i) = \sum_{k=1}^{n}p(n-k)A(k),
\end{equation}
where
\begin{equation}
A(n)=\sum_{d|n}a(d).
\end{equation}
\end{theorem}
\begin{proof}
Recognizing the right-hand side of the sum in Theorem \ref{Theorem1} as a Lambert series and $\frac{1}{(q;q)_\infty}$ as the generating function for $p(n)$, we obtain
\begin{equation}
\sum_{n=0}^{\infty}\sum_{\lambda \vdash n}\sum_{\lambda_i \in \lambda}a(\lambda_i)q^n = \frac{1}{(q;q)_{\infty}}\langle f\rangle_q  =\sum_{n=0}^{\infty}p(n)q^n \sum_{n=1}^{\infty}A(n)q^n,
\end{equation}
where $A(n)$ is as defined above and $f(\lambda)=\sum_{\lambda_i \in \lambda}a(\lambda_i)$. Taking a Cauchy product, reindexing, and comparing coefficients of $q^n$ yields the desired identity.
\end{proof}

\section{Sums over distinct parts}
\label{Sums2}

Now we address sums over the parts of partitions but without repetition. These sums enjoy similar relations to those given in Section \ref{Sums1}.
\begin{theorem}\label{Theorem3}
For $\lambda \in \mathcal P$ let $f_{\text{dist}}(\lambda)=\sum_{\substack{\lambda_i \in \lambda\\ \lambda_i\  \text{distinct}}} a(\lambda_i)$. Then, using the notation introduced above, we have
$$
\left< f_{\text{dist}} \right>_q=\sum_{n=1}^{\infty} a(n)q^n.
$$
\end{theorem}
\begin{proof}
Analogously to the proof of \ref{Theorem1}, we now study the product 
\begin{align*}
\phi_2(t,q) &:= (1+t^{a(1)} q+t^{a(1)}q^2+\ldots)(1+t^{a(2)} q^2+t^{a(2)}q^4+\ldots)\times\cdots \\
&\times(1+t^{a(n)} q^n+t^{a(n)}q^{2n}+\ldots)\times\cdots.
\end{align*}
The coefficient of each $q^n$ will again be a polynomial in $t$. However, exponent $m$ in $t^m$ will now be equal to $a(n_1)+a(n_2)+\ldots$ for distinct $n_1, n_2, \ldots, n_i$ in a partition of $n$ since we do not give different weights based on how many times each $n_i$ appears in a partition. Therefore, as before we have 
\begin{equation}
\frac{\partial \phi_2(t,q)}{\partial t}\Bigr|_{t=1} = \sum_{n=0}^{\infty}q^n \sum_{\lambda \vdash n}f_{\text{dist}}(\lambda) = \frac{1}{(q;q)_{\infty}}\langle f_{\text{dist}}\rangle_q.
\end{equation}
Evaluating the partial by summing each term in $\phi_2$ and taking the logarithmic derivative of  
\begin{equation}
\phi_2(t,q) = \prod_{n=1}^{\infty}\left(1+\frac{t^{a(n)}q^n}{1-q^n}\right)
\end{equation}
while noting $\phi_2(1,q) = \frac{1}{(q;q)_\infty}$ yields the desired theorem.
\end{proof}
The case for constant $a(n)$ was treated with a similar method in \cite[(Page 2)]{Hirschhorn1}.

\begin{theorem}\label{Theorem4}
We have the convolution identity
\begin{equation}
\sum_{\lambda \vdash n}\sum_{\substack{\lambda_i \in \lambda \\ \lambda_i \text{distinct} }}a(\lambda_i) = \sum_{k=1}^{n}p(n-k)a(n).
\end{equation}
\end{theorem}
\begin{proof}
Rewriting Theorem \ref{Theorem3} as
\begin{equation}
\sum_{n=0}^{\infty}\sum_{\lambda \vdash n}\sum_{\substack{\lambda_i \in \lambda \\ \lambda_i \text{distinct} }}a(\lambda_i)q^n =\sum_{n=0}^{\infty}p(n)q^n \sum_{n=1}^{\infty} a(n) q^n,
\end{equation}
taking a Cauchy product of the right-hand side, reindexing, and comparing coefficients of $q^n$, yields the desired identity.
\end{proof}

\section{Connections to multiplicative number theory}
\label{Stanley}
Here we draw parallels between the preceding theorems and objects in multiplicative number theory.
\begin{theorem}\label{MultCor}
Using the notation introduced above, we have the identity
$$
\sum_{\lambda \vdash n}\sum_{\substack{\lambda_i \in \lambda \\ \lambda_i \text{distinct} }}A(\lambda_i)=\sum_{\lambda \vdash n}\sum_{\lambda_i \in \lambda}a(\lambda_i) = \sum_{k=1}^{n}p(n-k)A(k),
$$
where
\begin{equation*}\label{3.5}
A(n)=\sum_{d|n}a(d).
\end{equation*}
In other words, the sum of $f$ over all parts in every partition of $n$ is equal to the divisor sum of $f$, evaluated over all distinct parts in every partition of $n$.
\end{theorem}
\begin{proof}
Looking at the convolution identities Theorem \ref{Theorem2} and Theorem \ref{Theorem4}, we see that pairs of functions that satisfy $A(n)=\sum_{d|n}a(n)$ lead to expressions of the form 
\begin{equation}
\sum_{\lambda \vdash n}\sum_{\substack{\lambda_i \in \lambda \\ \lambda_i \text{distinct} }}A(\lambda_i) = \sum_{k=1}^{n}p(n-k)\sum_{d|k}a(d) = \sum_{\lambda \vdash n}\sum_{\lambda_i \in \lambda}a(\lambda_i) .
\end{equation}
This completes the proof.
\end{proof}

This result generalizes Stanley's Theorem \cite[pp. 185-186]{Stanley}, which states that the total number of $1$s in all partitions of $n$ is equal to the sum of the number of distinct parts in each partition. Theorem \ref{MultCor} reduces to Stanley's Theorem if we let $a(n)$ be the indicator function for $1$, so that $A(n)=1$. 

Theorem \ref{MultCor} is quite interesting in that it links divisor sums, which are ubiquitous in multiplicative number theory, to sums over partitions of $n$, a structure central to additive number theory. For further reading at the intersection of these two branches of number theory see the works of K. Alladi, for instance \cite{Alladi1}. By specializing $a(n)$ in Theorem \ref{MultCor}, we obtain several apparently new corollaries.

\begin{corollary}
The number of partitions of $n$ having $1$ as a part is equal to $p(n-1)$. The number of squarefree parts with an even number of prime factors minus the number of squarefree parts with an odd number of prime factors, summed over every partition of $n$, is also given by $p(n-1)$.
\end{corollary}
\begin{proof}
Take $a(n)=\mu(n)$ in Theorem \ref{MultCor}, where $\mu(n)$ is the M{\"o}bius function, so that $A(n)$ is the indicator function for $1$ \cite[(Page 25)]{Apostol}. The convolution on the right-hand side then sums to $p(n-1)$, while $\sum_{\lambda \vdash n}\sum_{\substack{\lambda_i \in \lambda \\ \lambda_i \text{distinct} }}A(\lambda_i)$ then counts the number of $1$s in the distinct parts of each partition of $n$. This is equivalent to simply counting the number of partitions of $n$ which contain a $1$. We note that this can be directly proved by removing a $1$ from every partition of $n$ that contains it, since these will then form partitions of $n-1$ (George Andrews, personal communication, 2016).

The second statement follows from applying the definition of the M{\"o}bius function to $\sum_{\lambda \vdash n}\sum_{\lambda_i \in \lambda}\mu(\lambda_i)$. 
\end{proof}

\begin{corollary}We have the identity
\begin{equation}
\sum_{\lambda \vdash n}\sum_{\lambda_i \in \lambda}\lambda_i^\alpha=\sum_{\lambda \vdash n}\sum_{\substack{\lambda_i \in \lambda \\ \lambda_i \text{distinct} }} \sigma_\alpha(\lambda_i)=\sum_{k=1}^{n}\sigma_{\alpha}(k)p(n-k).
\end{equation}
\end{corollary}
\begin{proof}
Take $a(n)=n^\alpha$ in Theorem \ref{MultCor}.
\end{proof}

This forms a natural generalization of the classical identity, originally derived by Euler \cite[(Page 323)]{Apostol},
\begin{equation}\label{4.5}
np(n)=\sum_{k=1}^{n}\sigma_{1}(k)p(n-k),
\end{equation}
since $\sum_{\lambda \vdash n}\sum_{\lambda_i \in \lambda}\lambda_i = \sum_{\lambda \vdash n}n = n p(n)$. Both the inclusion of the parameter $\alpha$ and the connection to sums over distinct parts are apparently new.

\begin{corollary}We have the identity
\begin{equation}
\sum_{\lambda \vdash n}\sum_{\lambda_i \in \lambda}J_\alpha(\lambda_i)=\sum_{\lambda \vdash n}\sum_{\substack{\lambda_i \in \lambda \\ \lambda_i \text{distinct} }} \lambda_i^\alpha=\sum_{k=1}^{n}J_{\alpha}(k)p(n-k).
\end{equation}
\end{corollary}
\begin{proof}
Take $a(n)=J_\alpha(n)$ in Theorem \ref{MultCor}, where $J_\alpha(n)$ is the Jordan totient function \cite[(48)]{Apostol}, so that $A(n)=n^\alpha$.
\end{proof}

The forms a natural dual of Euler's convolution, since the $\sigma_\alpha$ convolution arises from considering $\alpha$th moments of all parts of the partitions of $n$, while the $J_\alpha$ convolution arises from considering $\alpha$th moments of all distinct parts of the partitions of $n$. We also note that taking $\alpha=1$ provides an analog of \eqref{4.5} for the Euler totient function.

\begin{corollary}
Let $Q(n)$ denote the number of squarefree parts in all partitions of $n$. Then $Q(n) = \sum_{k=1}^{n}p(n-k)2^{\omega(k)}$, where $\omega(n)$ gives the number of distinct prime factors of $n$. Therefore, $Q(n) \equiv p(n-1) \pmod 2$.
\end{corollary}
\begin{proof}
Take $a(n)=\mu^2(n)$ in Theorem \ref{MultCor}, so that $A(n)=2^{\omega(n)}$  \cite[(Page 45)]{Apostol}, where $\omega(n)$ is as given above. We have that $a(n)$ is the indicator function for squarefree numbers, so $\sum_{\lambda \vdash n}\sum_{\lambda_i \in \lambda}\mu^2(n)$ gives the number of squarefree parts in all partitions of $n$. The congruence arises from noting that $2^{\omega(k)} \equiv 0 {\pmod 2}$ unless $k=1$.
\end{proof}

\begin{corollary}
The sum of the number of distinct squares in every partition of $n$ is given by
\begin{equation}
\sum_{k=1}^{\left \lfloor{\sqrt{n}}\right \rfloor}p(n-k^2).
\end{equation}
\end{corollary}
\begin{proof}
Take $a(n)=\lambda(n)$ in Theorem \ref{MultCor}, where $\lambda(n):=(-1)^{\Omega(n)}$ is Liouville's function \cite[pp. 37-38]{Apostol} (note we are using $\lambda$ in two different ways now). Then $A(n)$ is the indicator function for the squares \cite[(Page 38)]{Apostol}. This implies $\sum_{k\leq n, k\text{ square}}p(n-k)$ is the sum of the number of distinct squares in every partition of $n$. Rewriting the summation completes the proof.
\end{proof}

\begin{corollary}\label{Cor10}We have the identity
\begin{equation}\label{3.7}
\sum_{\lambda \vdash n}\sum_{\lambda_i \in \lambda}\Lambda(\lambda_i)=\sum_{\lambda \vdash n}\sum_{\substack{\lambda_i \in \lambda \\ \lambda_i \text{distinct} }} \log(\lambda_i)=\sum_{k=1}^{n}p(n-k) \log{k}.
\end{equation}
Furthermore, $\prod_{k=1}^{n}k^{p(n-k)}$ equals the product of all distinct parts in every partition of $n$.
\end{corollary}
\begin{proof}
Take $a(n)=\Lambda(n)$ in Theorem \ref{MultCor}, where $\Lambda(n)$ is the von Mangoldt function, which equals $\operatorname{log} p$ if $n$ is an integer power of the prime $p$, and is 0 otherwise. Then, $A(n)=\log n$ \cite[(Page 32)]{Apostol}, which yields the corollary. By exponentiating the sum over distinct $\lambda_i$ and the right-hand side, we obtain the second part of the corollary.
\end{proof}

An open question for future research is now to provide combinatorial or alternative proofs for these formulas, or apply Theorem \ref{MultCor} to any other divisor sum and find new consequences. For instance, studying the prime factorization of $\prod_{k=1}^{n}k^{p(n-k)}$ or the role of the von Mangoldt function in Corollary \ref{Cor10} could yield information about the multiplicities of primes in the parts of the partitions of $n$.

\begin{acknowledgements}
I would like to thank Christophe Vignat for guiding me through my first forays into mathematics, and Armin Straub for first introducing me to Stanley's Theorem and the current literature on partitions. I would also like to thank George Andrews, Robert Schneider, and the anonymous referee for helpful comments. Last but not least, I'd like to thank Ellicott 3 for always providing me with (questionable) inspiration.
\end{acknowledgements}



\end{document}